\documentclass[a4paper,conference,top=2cm]{IEEEtran}
\IEEEoverridecommandlockouts
\usepackage{etoolbox}
\makeatletter
\patchcmd{\@makecaption}
  {\scshape}
  {}
  {}
  {}
\makeatother

\usepackage{caption}
\def\BibTeX{{\rm B\kern-.05em{\sc i\kern-.025em b}\kern-.08em
    T\kern-.1667em\lower.7ex\hbox{E}\kern-.125emX}}

\usepackage{cite}
\usepackage{amsmath,amssymb,amsfonts}
\usepackage{algorithmic}
\usepackage{graphicx}
\usepackage{textcomp}
\usepackage{xcolor}

\usepackage{stfloats}
\def\BibTeX{{\rm B\kern-.05em{\sc i\kern-.025em b}\kern-.08em
    T\kern-.1667em\lower.7ex\hbox{E}\kern-.125emX}}

\usepackage{import}

\usepackage[utf8]{luainputenc}

\usepackage[T1]{fontenc}

\usepackage{graphicx}

\usepackage[english]{babel}
\addto\captionsenglish{\renewcommand{\figurename}{Fig.}}
\addto\captionsenglish{\renewcommand{\tablename}{Tab.}}
\usepackage{csquotes}

\usepackage{newtxtext}
\usepackage{amsthm}
\usepackage[slantedGreek]{newtxmath}
\usepackage[OMLmathsfit]{isomath}
\DeclareMathAlphabet{\mathbfsf}{\encodingdefault}{\sfdefault}{bx}{n}
\usepackage{bm}
\usepackage{envmath}
\usepackage{mathtools}
\usepackage{commath}
\usepackage{siunitx}

\usepackage[caption=false,font=footnotesize]{subfig}

\usepackage{booktabs}
\usepackage{footmisc}  

\usepackage{url}

\theoremstyle{definition}

\theoremstyle{plain}

\theoremstyle{remark}

\usepackage{lineno}
\modulolinenumbers[5]
\usepackage{umoline}

\usepackage{pgfplots}
\usepackage{pgfplotstable}
\pgfplotsset{compat=newest}
\pgfplotsset{plot coordinates/math parser=false}
\newlength\figureheight
\newlength\figurewidth
\pgfplotsset{every axis plot/.append style={line width=1.5pt},
    legend style={font=\footnotesize, 
        text height=1.0ex,
        draw=black,
        fill=white,
        legend cell align=left}}

\usepackage{hyperref} 
\usepackage[english]{cleveref}

\Crefname{defn}{definition}{definitions}
\Crefname{defn}{Definition}{Definitions}

\Crefname{asm}{assumption}{assumptions}
\Crefname{asm}{Assumption}{Assumptions}

\crefname{lem}{lemma}{lemmas} 
\Crefname{lem}{Lemma}{Lemmas}

\crefname{prop}{proposition}{propositions} 
\Crefname{prop}{Proposition}{Propositions}

\crefname{thm}{theorem}{theorms} 
\Crefname{thm}{Theorem}{Theorms}

\crefname{cor}{corollary}{corollaries}
\Crefname{cor}{Corollary}{Corollaries}
\newcounter{subequation}
\newlength\mtabskip\mtabskip=-1.25cm

\def\mtabLong{long}

\newcommand{\mr}{\mathrm}

\newcommand{\veg}[1]{\bm{#1}}     
\newcommand{\mat}[1]{\mathsfbfit{#1}} 
\renewcommand{\vec}[1]{\mathsfbfit{#1}} 
\newcommand{\vecop}[1]{\bm{\mathcal{#1}}} 

\newcommand{\dd}{\mathrm{d}}  

\newcommand{\T}{\mr{T}}

\newcommand\restr[2]{{
        \left.\kern-\nulldelimiterspace 
        #1 
        \vphantom{|} 
        \right|_{#2} 
}}

\newcommand\rst[3]{{
        \left.\kern-\nulldelimiterspace 
        #1 
        \vphantom{|} 
        \right|_{#2}^{#3} 
}}

\usepackage{acro}

\DeclareAcronym{DG}
{
    short = DG ,
    long = discontinuous Galerkin
}

\DeclareAcronym{ACA}
{
    short = ACA ,
    long = adaptive cross approximation
}

\DeclareAcronym{EFIE}
{
    short =  EFIE ,
    long = electric field integral equation
}

\DeclareAcronym{MFIE}
{
    short =  MFIE ,
    long = magnetic field integral equation
}

\DeclareAcronym{CFIE}
{
    short =  CFIE ,
    long = combined field integral equation
}

\DeclareAcronym{MUIE}
{
    short =  MUIE ,
    long = Müller integral equation
}

\DeclareAcronym{PMCHWT}
{
    short =  PMCHWT ,
    long = Poggio-Miller-Chang-Harrington-Wu-Tsai integral equation
}

\DeclareAcronym{SPD}
{
    short =  SPD ,
    long = {symmetric, positive definite}
}

\DeclareAcronym{SPSD}
{
    short =  SPD ,
    long = {symmetric, positive semi-definite}
}

\DeclareAcronym{PEC}
{
    short =  PEC ,
    long = perfectly electrically conducting
}

\DeclareAcronym{RWG}
{
    short = RWG ,
    long = Rao-Wilton-Glisson
} 

\DeclareAcronym{BC}
{
    short = BC ,
    long = Buffa-Christiansen
}

\DeclareAcronym{SVD}
{
    short = SVD ,
    long = singular value decomposition
}

\DeclareAcronym{CG}
{
    short = CG ,
    long = conjugate gradient
} 

\DeclareAcronym{PCG}
{
    short = PCG ,
    long = preconditioned conjugate gradient
} 

\DeclareAcronym{CGS}
{
    short = CGS ,
    long = conjugate gradient squared
}

\DeclareAcronym{CMP}
{
    short = CMP ,
    long = Calderón multiplicative preconditioner
} 

\DeclareAcronym{RFCMP}
{
    short = RF-CMP ,
    long = refinement-free Calderón multiplicative preconditioner
} 

\DeclareAcronym{HPD}
{
    short = HPD ,
    long = {Hermitian, positive definite}
} 

\DeclareAcronym{RHS}
{
    short = RHS ,
    long = {right-hand side}
}

\DeclareAcronym{PW}
{
    short = PW ,
    long = {plane wave}
} 

\DeclareAcronym{HD}
{
    short = HD ,
    long = {Hertzian dipole}
} 

\DeclareAcronym{FF}
{
    short = FF ,
    long = {far-field}
} 

\DeclareAcronym{NF}
{
    short = NF ,
    long = {near-field}
}  

\newcolumntype {n}{c}
\newcolumntype {N}{>{\small}c}
\newcolumntype {L}{>{\small}l}
\newcolumntype {F}{>{\footnotesize}c}
\newcolumntype {v}[1]{>{\raggedright \hspace {0pt}} p {#1}}
\newcolumntype {V}[1]{>{\small \raggedright \hspace {0pt}} p {#1}}
\newcolumntype{d}[1]{>{\DC@{.}{.}{#1}}c<{\DC@end}}

\newcolumntype{R}[1]{%
    >{\begin{turn}{90}\begin{minipage}{#1}\small\raggedright\hspace{0pt}}l%
            <{\end{minipage}\end{turn}}%
}

\NewDocumentCommand{\TA}{o}{
    \IfNoValueTF {#1} {%
        \vecop T_{\kern-2pt\mr{A}}
    }
    {
        \vecop T_{\kern-2pt\mr{A},#1}
    }
}

\NewDocumentCommand{\TPhi}{o}{
    \IfNoValueTF {#1} {%
        \vecop T_{\kern-2pt\Phiup}
    }
    {
        \vecop T_{\kern-2pt\Phiup,#1}
    }
}

\NewDocumentCommand{\matTA}{o}{
    \IfNoValueTF {#1} {%
        \mat T_\mr{A}   
        }
    {
        \mat T_{\mr{A},#1}
    }
}

\NewDocumentCommand{\matTPhi}{o}{
    \IfNoValueTF {#1} {%
        \mat T_\Phiup   
        }
    {
        \mat T_{\Phiup,#1}
    }
}

\NewDocumentCommand{\MSL}{o}{
    \IfNoValueTF {#1} {%
        \veg \Psi_\mr{SL}
        }
    {
        \veg \Psi_{\mr{SL},#1}
    }
}

\NewDocumentCommand{\MDL}{o}{
    \IfNoValueTF {#1} {%
        \veg \Psi_\mr{DL}
        }
    {
        \veg \Psi_{\mr{DL},#1}
    }
}

\NewDocumentCommand{\PA}{o}{
    \IfNoValueTF {#1} {%
        \veg \Psi_\mr{A}
        }
    {
        \veg \Psi_{\mr{A},#1}
    }
}

\NewDocumentCommand{\PPhi}{o}{
    \IfNoValueTF {#1} {%
        \veg \Psi_{\Phiup}
        }
    {
        \veg \Psi_{\Phiup,#1}
    }
}

\usepackage{todonotes}
\usepackage{comment}

\makeatletter

\newcounter{authr}
\newcommand{\authr}[2][]{
   \stepcounter{authr}
   \@namedef{authr@\theauthr}{#2}
   \@namedef{authrlabel@\theauthr}{#1}
}

\newcounter{address}
\newcommand{\address}[2][]{
   \stepcounter{address}
   \@namedef{address@\theaddress}{#2}
   \@namedef{addresslabel@\theaddress}{#1}
}

\newcommand{\alsep}{and}

\def\newmaketitle{\par%
  \begingroup%
  \normalfont%
  \def\thefootnote{}
  \def\footnotemark{}
  \let\@makefnmark\relax
  \footnotesize
  \footnotesep 0.7\baselineskip
  \normalsize%
  \twocolumn[\thenewmaketitle\@IEEEaftertitletext]%
  \if@IEEEusingpubid
     \enlargethispage{-\@IEEEpubidpullup}%
  \fi
  \endgroup
  \setcounter{footnote}{0}\let\maketitle\relax\let\@maketitle\relax
  \gdef\@thanks{}%
  \let\thanks\relax}

\def\thenewmaketitle{
  \newpage
  \begin{center}%
    \vskip0.2em{\Huge\@IEEEcompsoconly{\sffamily}\@IEEEcompsocconfonly{\normalfont\normalsize\vskip 2\@IEEEnormalsizeunitybaselineskip
   \bfseries\large}\@title\par}\vskip1.0em\par%
    \vspace{1ex}
    \newcounter{c@authr}
    \newcounter{c@tmp}
    \ifthenelse{\value{authr}=2}{%
      \newcommand{\liand}{ and }}{%
      \newcommand{\liand}{, and }}
    \ifthenelse{\value{address}<2}{%
      \@nameuse{authr@1}%
      \stepcounter{c@authr}%
      \whiledo{\value{c@authr}<\value{authr}}{%
        \setcounter{c@tmp}{\value{authr}}%
        \addtocounter{c@tmp}{-\value{c@authr}}%
        \ifthenelse{\value{c@tmp}=1}{%
          \renewcommand{\alsep}{\liand}}{\renewcommand{\alsep}{, }}%
        \stepcounter{c@authr}\alsep \@nameuse{authr@\thec@authr}}\\%
    }
    {
      \@nameuse{authr@1}${}^{(\ref{\@nameuse{authrlabel@1}})}$%
      \stepcounter{c@authr}%
      \whiledo{\value{c@authr}<\value{authr}}{%
      \setcounter{c@tmp}{\value{authr}}%
      \addtocounter{c@tmp}{-\value{c@authr}}%
      \ifthenelse{\value{c@tmp}=1}{%
        \renewcommand{\alsep}{\liand}}{\renewcommand{\alsep}{, }}%
      \stepcounter{c@authr}\alsep \@nameuse{authr@\thec@authr}%
        ${}^{(\ref{\@nameuse{authrlabel@\thec@authr}})}$%
      }
    }
    \vspace{0.2ex}

    \ifthenelse{\value{address}>0}{%
      \ifthenelse{\value{address}=1}{
        {\@nameuse{address@1}}
      }
      {
        \newcounter{c@address}

        \begin{center}
        \whiledo{\value{c@address}<\value{address}}
        {
          \refstepcounter{c@address}
            ${}^{(\thec@address)}$\,%
              \label{\@nameuse{addresslabel@\thec@address}}%
              \@nameuse{address@\thec@address}\\ %
        }
        \end{center}
      } 
    }
    {
      \relax
    }
  \end{center}
}

\makeatother

\title{Linear-in-Complexity Computational Strategies for Modeling and Dosimetry at TeraHertz
}

\authr[org1]{Viviana Giunzioni}
\authr[org1]{Giuseppe Ciacco}
\authr[org2]{Clément Henry}
\authr[org2]{Adrien Merlini}
\authr[org1]{Francesco P. Andriulli}

\address[org1]{Department of Electronics and Telecommunications, Politecnico di Torino, Italy}
\address[org2]{Microwaves Department, IMT Atlantique, Brest, France}

\begin{document}

\newmaketitle

\begin{abstract}
This work presents a fast direct solver strategy allowing full-wave modeling and dosimetry at terahertz (THz) frequencies. The novel scheme leverages a preconditioned combined field integral equation together with a regularizer for its elliptic spectrum to enable its compression into a non-hierarchical skeleton, invertible in quasi-linear complexity. Numerical results will show the effectiveness of the new scheme in a realistic skin modeling scenario. 
\end{abstract}


\section{Introduction}
With the technological advances in THz technology, a growing number of interdisciplinary applications in the THz range have emerged and gained popularity in the last two decades within areas ranging from security, to communications, or biomedicine \cite{chen2022terahertz}. As the impact of THz devices in our societies grows, accurately assessing the effects of THz waves on the human body gains crucial importance \cite{cherkasova2020effects}. Hence the need for exposure analyses that aim at quantifying the amount of energy absorbed by biological tissues subject to electromagnetic radiations \cite{alekseev2008millimeter}.

Preliminary dosimetry assessments are a fundamental phase during the design of THz equipments, to guarantee their compliance with the limits on the power absorbed by human tissues set by international agencies \cite{icnirp2020guidelines}.
Exposure measurements are often challenging to perform, especially in the near field, but this challenge can be, in part, sidestepped by reliable and accurate numerical dosimetric assessments, when they are within reach. However numerical modeling at THz also comes with its own set of complications.


On the one hand, many of the solvers proposed in the literature employ approximations of the Maxwell's system, suitable to the high frequency regime considered, or apply geometrical simplification to make use of proper analytic solutions. However, application of these approximations can degrade the solution accuracy, and potentially compromise the reliability of the dosimetric analyses.
On the other hand, full-wave models leverage the original Maxwell system and can be applied to arbitrarily complex geometries, 
but the  higher computational costs incurred can become prohibitive. In addition, they suffer from numerical issues, such as ill-conditioning or spurious resonances at high frequencies \cite{adrian2021electromagnetic} that need to be handled to obtain reliable results.

We propose here a novel full-wave approach, well suited to modeling reflection and absorption of THz waves by biological samples. Being a fast direct solution strategy, this approach allows for the efficient solution of the THz problems for multiple exposures at once, with a complexity which grows only quasi-linearly with the number of unknowns, that is, with increasing frequency. This is obtained by first defining a proper set of boundary integral equations and leveraging a tailored preconditioning scheme, resulting in a well-conditioned system of linear equations freed from spurious resonances. This formulation is then coupled with a recently proposed fast inversion strategy \cite{consoli2022fast}, that relies on the compression of the elliptic spectrum of the boundary operator into a rank-deficient skeleton form and on the use of the Woodbury matrix identity \cite{henderson1981deriving}.


\section{Background and Notation}
Dosimetry analyzes aim at assessing the amount of energy absorbed by the human body when exposed to an electromagnetic radiation. This estimation can be performed by numerically simulating the response of the biological tissue to the impinging field through a full-wave electromagnetic solver.
In this work we employ the two-dimensional approximation, that assumes tha invariance of the geometries and field along an axis $\hat{\veg z}$.
This lends itself well to the case under study given the large dimensions of some body parts compared to THz wavelengths. However, this approximation is not suited to modeling all body parts.

Based on the representation theorem \cite{nedelec2001acoustic}, different boundary integral equations (BIEs) can be set up to numerically model the time-harmonic electromagnetic scattering and absorption of a penetrable body.
Given a two-dimensional domain $\Omega$ with boundary $\Gamma\coloneqq\uppartial\Omega$ characterized by the outgoing normal field $\hat{\veg n}$, the boundary integral operators \cite{nedelec2001acoustic}
\begin{align}
\left(\mathcal{S}_k^\Gamma\psi\right)(\veg{r}) &\coloneqq  \int_{\Gamma}G_k(\veg{r}- \bm{r}') \psi(\veg{r}') \dd S(\veg{r}') \,,
\label{eqn:Sop}\\
\left(\mathcal{D}_k^\Gamma\psi\right)(\veg{r}) &\coloneqq  \mathrm{p.v.} \int_{\Gamma}\uppartial_{\veg{n}'} G_k(\veg{r}- \veg{r}') \psi(\veg{r}') \dd S(\veg{r}') \,,
\label{eqn:Dop}\\
\left(\mathcal{D}^{* \Gamma}_k\psi\right)(\veg{r}) &\coloneqq  \mathrm{p.v.} \int_{\Gamma}\uppartial_{\veg{n}} G_k(\veg{r}- \veg{r}') \psi(\veg{r}') \dd S(\veg{r}') \,,
\label{eqn:Dsop}\\
\left(\mathcal{N}_k^\Gamma\psi\right)(\veg{r}) &\coloneqq  -\mathrm{f.p.} \int_{\Gamma}\uppartial_{\veg{n}}\uppartial_{\veg{n}'}G_k(\veg{r}- \veg{r}') \psi(\veg{r}')\dd S(\veg{r}')\,,
\label{eqn:Nop}
\end{align}
which are respectively the single layer, double layer, adjoint double layer, and hypersingular operator, constitute the building blocks of any of these formulations. The notations $\text{p.v.}$ and $\text{f.p.}$ indicate the Cauchy principal value and the Hadamard finite part.
We denote by $G_k$ the two-dimensional Green's function in free-space
\begin{equation}
    G_k(\veg{r}- \veg{r}') = -\frac{j}{4}H_0^{(2)}\left( k ||\veg{r}- \veg{r}'|| \right)\,,
\end{equation}
where $H_0^{(2)}$ is the Hankel function of the second kind with order zero as defined in \cite{olver2010nist}.

Moreover, numerical exposure assessments also require the \textit{a priori} definition of a realistic model of the tissue under study, both in terms of geometry and dielectric permittivity. 
Research on THz external dosimetry is often focused on the skin \cite{alekseev2008millimeter,haider2022highresolution}, as THz impinging field is absorbed by this organ.
Different geometrical models of the skin have been proposed \cite{wang2021thz} to accurately reproduce the human anatomy. They usually aim at modeling the stratification of compartments with different physical properties, such as the stratum corneum, the epidermis, and the dermis layers, sometimes even modeling anisotropies and depth-varying water percentage \cite{wang2021thz}, at the cost of increasing model complexity. For the sake of simplicity, in this work we employ a single-dielectric model. Following the double Debye model \cite{kindt1996farinfrared}, the permittivity of the skin as a function of the frequency is modeled as
\begin{equation}
    \epsilon_r(\omega) = \epsilon_\infty+\frac{\epsilon_s-\epsilon_2}{1+j\omega\tau_1}+\frac{\epsilon_2-\epsilon_\infty}{1+j\omega\tau_2}\,,
    \label{eqn:debye}
\end{equation}
with parameters $\epsilon_\infty=$ \num{3}, $\epsilon_s=$ \num{60}, $\epsilon_2=$ \num{3.6}, $\tau_1=$ \SI{10}{ps}, and $\tau_2=$ \SI{0.2}{ps} \cite{pickwell2004simulation}. The validity of this approximation has been demonstrated in previous works \cite{pickwell2004vivo}, which however have also highlighted a limitation of the model when applied to dry skin and, in general, to tissues characterized by low water contents.

\section{Fast Direct Solver Strategy for THz Dosimetry}
We propose here a fast direct solver strategy for modeling the electromagnetic response of a biological tissue of boundary $\Gamma_s$ to an excitation realized by means of a metallic body of boundary $\Gamma_m$. It is based on a composite formulation made up of the combined field integral equation (CFIE) for perfect electric conductor (PEC) materials \cite{mavtz1978hfield} and of the Poggio-Miller-Chang-Harrington-Wu-Tsai (PMCHWT) equation for penetrable media \cite{poggio1973integral}.
As is sometimes done in the literature, we assume that the coupling terms between the metallic and the
dielectric objects can be neglected \cite{ziane2020antenna,sacco2021antenna}.
In the case where both objects (i.e., the metallic and the dielectric ones) are subject to a TM polarized field, the resulting system of integral equations is given in \cref{eqn:integraleq} at the bottom of the page. Similar results can be found for different polarizations. 

\begin{figure*}[b]
\normalsize
\setcounter{equation}{6}
\begin{equation}
\begin{cases}
    \mathcal{S}_{k_0}^{\Gamma_m}  j_{z,m}(\veg r) + \left(\frac{1}{2}\mathcal{I}+\mathcal{D}^{* \Gamma_m}_{k_0} \right)j_{z,m}(\veg r) = \frac{1}{jk_0\eta} E_z^{\text{inc},m}(\veg r) + H_t^{\text{inc},m}(\veg r) , & \text{$\veg r \in \Gamma_m$}\\
    \left(-jk_0\eta_0\mathcal{S}_{k_0}^{\Gamma_s} - jk_1\eta_1\mathcal{S}_{k_1}^{\Gamma_s} \right) j_{z,s}(\veg r) + \left( \mathcal{D}_{k_0}^{\Gamma_s}+\mathcal{D}_{k_1}^{\Gamma_s}\right) m_{t,s}(\veg r) = E_z^{\text{inc},s}(\veg r) , & \text{$\veg r \in \Gamma_s$}\\
    -\left( \mathcal{D}_{k_0}^{*\Gamma_s}+\mathcal{D}_{k_1}^{*\Gamma_s}\right) j_{z,s}(\veg r) + \left(-1/(jk_0\eta_0)\mathcal{N}_{k_0}^{\Gamma_s} - 1/(jk_1\eta_1)\mathcal{N}_{k_1}^{\Gamma_s} \right) m_{t,s}(\veg r) = H_t^{\text{inc},s}(\veg r) , & \text{$\veg r \in \Gamma_s$}
\end{cases}
\label{eqn:integraleq}
\end{equation}
\vspace*{4pt}
\end{figure*}

In these equations, we denote by the subscript $_0$ the quantities related to the exterior medium, which can be assumed to be the air, and by the subscript $_1$ the ones related to the interior, penetrable, medium. The exterior and interior wavenumbers are denoted as $k_0=\omega\sqrt{\epsilon_0\mu_0}$ and $k_1=\omega\sqrt{\epsilon_1\mu_1}$, while $\eta_{0/1} = \sqrt{\mu_{0/1}/\epsilon_{0/1}}$ are the characteristic impedances of the exterior or interior medium. $(E^{\text{inc},m},H^{\text{inc},m})$ and $(E^{\text{inc},s},H^{\text{inc},s})$ are the electromagnetic fields incident over $\Gamma_m$ and $\Gamma_s$ respectively, further separated into the transversal and longitudinal components, denoted by $_t$ and $_z$. In the following, we will denote by the subscripts $_m$ and $_s$ quantities related to $\Gamma_m$ and $\Gamma_s$ respectively.

The unknowns in \cref{eqn:integraleq} are the surface equivalent currents defined on the metallic and dielectric boundaries. They are of electric type only in the former case, $j_{z,m}$, and of both electric and magnetic type in the latter, $j_{z,s}$ and $m_{t,s}$.
In particular, by superimposing the radiation provided by these currents, it is possible to retrieve the scattered electromagnetic field, to be summed to the incident field in order to determine the resulting electric and magnetic fields everywhere and, in particular, inside the biological sample.

The application of a Galerkin discretization scheme, based on the approximation of the unknown currents as linear combinations of $N_{m/s}$ piecewise linear basis functions $\lambda_i(\veg r)$ defined on a mesh of the boundary $\Gamma_{m/s}$, as $j_{z,m/s} \simeq \sum_{i=1}^{N_{m/s}} (\vec j_{z,m/s})_i \lambda_i$ and $m_{t,s} \simeq \sum_{i=1}^{N_{m/s}} (\vec m_{t,s})_i \lambda_i$, results in the linear system of equations
\begin{equation}
    \begin{pmatrix}
        \mat C & \vec 0 & \vec 0 \\
        \vec 0 & \mat P_{11} & \mat P_{12} \\
        \vec 0 & \mat P_{21} & \mat P_{22}
    \end{pmatrix}
    \begin{pmatrix}
        \vec j_{z,m} \\
        \vec j_{z,s} \\
        \vec m_{t,s}
    \end{pmatrix} = 
    \begin{pmatrix}
        \vec e_{z,m}/(jk_0\eta_0) + \vec h_{t,m}\\
        \vec e_{z,s} \\
        \vec h_{t,s}
    \end{pmatrix}\,.
    \label{eqn:lineareq}
\end{equation}
In the above system, 
\begin{align}
    (\vec e_{z,m/s})_i&=\left( \lambda_i, E_z^{\text{inc},m/s} \right)_{L^2(\Gamma_{s/m})}\quad\\
    (\vec h_{t,m/s})_i&=\left( \lambda_i, H_t^{\text{inc},m/s} \right)_{L^2(\Gamma_{s/m})}\,,
\end{align}
and the matrices $\mat C$, $\mat P_{11}$, $\mat P_{12}$, $\mat P_{21}$, and $\mat P_{22}$ are defined as
\begin{align}
    \mat C &= \mat{S}_{k_0}^{\Gamma_m}  + \frac{1}{2}\mat{G}^{\Gamma_m}+\mat{D}^{* \Gamma_m}_{k_0} \\
    \mat P_{11} &= -jk_0\eta_0\mat{S}_{k_0}^{\Gamma_s} - jk_1\eta_1\mat{S}_{k_1}^{\Gamma_s} \\
    \mat P_{12} &= \mat{D}_{k_0}^{\Gamma_s}+\mat{D}_{k_1}^{\Gamma_s} \\
    \mat P_{21} &= -\left(\mat{D}_{k_0}^{*\Gamma_s}+\mat{D}_{k_1}^{*\Gamma_s}\right) \\
    \mat P_{22} &= -1/(jk_0\eta_0)\mat{N}_{k_0}^{\Gamma_s} - 1/(jk_1\eta_1)\mat{N}_{k_1}^{\Gamma_s}\,,
\end{align}
where we have used the generic notation $(\mat O_{k}^\Gamma)_ij = \left(\lambda_i,  \mathcal{O}_{k}^\Gamma \lambda_j\right)_{L^2(\Gamma)}$, where $\mathcal{O}$ stands for one of $\{\mathcal{S}, \mathcal{D}, \mathcal{D}^*, \mathcal{N}\}$. The gram matrix $\mat G^\Gamma$ is obtained as $(\mat G^\Gamma)_{ij} = (\lambda_i,\lambda_j)_{L^2(\Gamma)}$.

As a consequence of the fact that the metallic radiator is electrically much larger than the biological sample under study, a significantly higher number of basis functions is required for the discretization of the unknown currents on its boundary $\Gamma_m$ (following the Nyquist sampling principle). Hence, we infer that the numerical solution of the linear system resulting from the discretization of the CFIE is the bottleneck, in terms of time and memory required, towards the solution of the entire system \eqref{eqn:lineareq}, both directly or iteratively. To alleviate this computational burden, we propose here to extend the Calder\'{o}n preconditioned scheme presented in \cite{andriulli2015high,consoli2022fast} for $\mat C$  and to extend the fast direct solver tailored for the resulting well-conditioned operator recently proposed in \cite{consoli2022fast}. In particular, we define the Calder\'{o}n stabilized operator matrix as
\begin{align}
    \mat C_p \coloneqq \,\,&\mat{N}_{\tilde{k}_0}^{\Gamma_m}\left(\mat{G}^{\Gamma_m}\right)^{-1}\mat{S}_{k_0}^{\Gamma_m}  + \nonumber\\
    &\left(\frac{1}{2}\mat{G}^{\Gamma_m}-\mat{D}^{* \Gamma_m}_{\tilde{k}_0}\right)\left(\mat{G}^{\Gamma_m}\right)^{-1}\left(\frac{1}{2}\mat{G}^{\Gamma_m}+\mat{D}^{* \Gamma_m}_{k_0}\right)\,
\end{align}
where, following the approach introduced in \cite{darbas2006generalized}, $\tilde{k}_0 \coloneqq k_0 - j 0.4 k_0^{1/3} a^{-2/3}$, with $a$ evaluated as a suitable average of the radius of curvature along $\Gamma_m$.
Then, following the procedure described in \cite{consoli2022fast}, we express $\mat C_p$ as the sum $\mat C_p = \mat C_{p,\text{c}}+\mat C_{p,\text{ext}}$, where $\mat C_{p,\text{c}}$ is the circular counterpart of $\mat C_p$ discretized over an equi-perimeter circular boundary. We employ at this point an adaptive randomized algorithm, such as the one presented in \cite{halko2011finding}, to compute a skeleton form of $\mat C_{p,\text{ext}}$ as
\begin{equation}
    \mat C_{p,\text{ext}} = \mat C_p-\mat C_{p,\text{c}}\simeq \mat U \mat V^\T\,.
\end{equation}
Given the spectral properties of matrix $\mat C_{p,\text{ext}}$, the rank of the skeleton $\mat U \mat V^\T$ grows only approximately as $k_0^{1/3}$ toward the high frequency limit. As a consequence, by applying a proper acceleration technique such as the fast multiple method (FMM) \cite{coifman1993fast}, the solution of the system, for any number of right hand sides, can be obtained efficiently, in quasi-linear complexity, by directly evaluating the inverse \cite{henderson1981deriving}
\begin{equation}
    \mat C_p^{-1} = \mat C_{p,\text{c}}^{-1} - \mat C_{p,\text{c}}^{-1} \mat U \left(\mat I + \mat V^\T \mat C_{p,\text{c}}^{-1}\mat U \right)^{-1}\mat V^\T \mat C_{p,\text{c}}^{-1}\,.
    \label{eqn:inverse}
\end{equation}
In particular, after noticing that all operations involving circulant matrices are computed rapidly via the use of the fast Fourier transform (FFT) algorithm, we recognize that the complexity of evaluating \eqref{eqn:inverse} scales in frequency approximately as $k_0^{4/3}$, with an overhead complexity with respect to the linear one determined by the skeleton rank increase.

\section{Numerical results}
\begin{figure}
\centerline{\includegraphics[width=1\columnwidth]{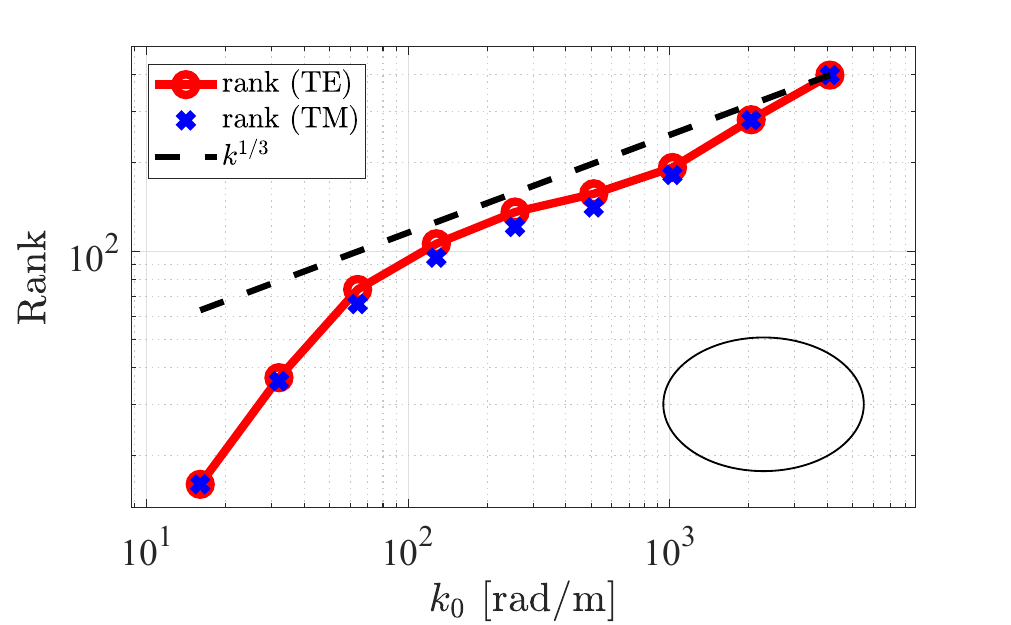}}
\caption{Rank of the skeleton form $\mat U \mat V^\T$ evaluated over an ellipse with aspect ratio \num{1.5} and perimeter \num{2}$\pi$ \unit{m} as a function of the free-space wavenumber $k_0$.}
\label{fig:ellipse_rank}
\end{figure}

\begin{figure}
\centerline{\includegraphics[width=1\columnwidth]{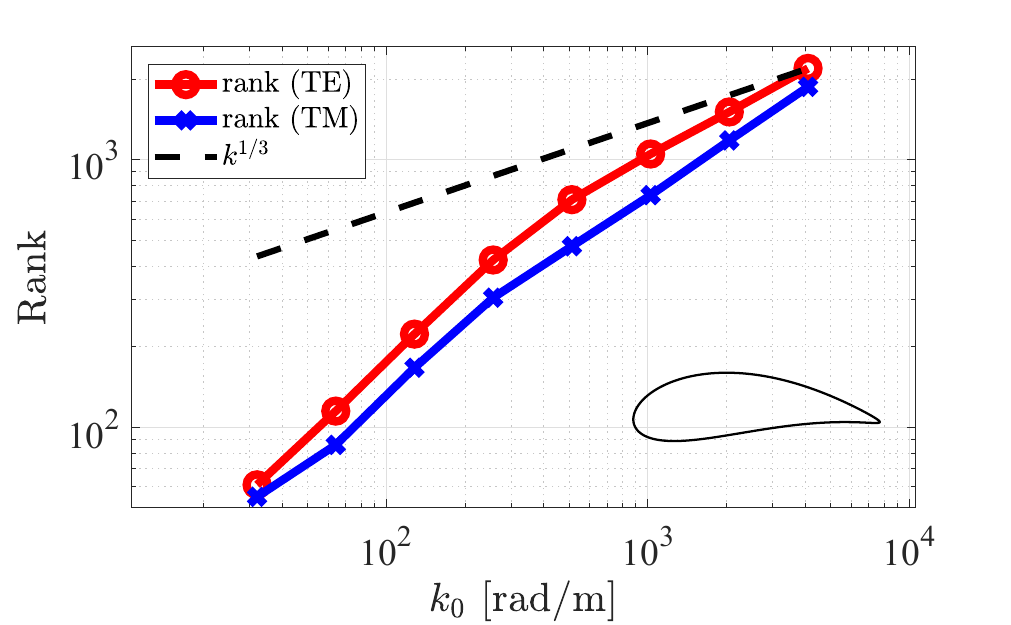}}
\caption{Rank of the skeleton form $\mat U \mat V^\T$ evaluated over an airfoil with perimeter \num{2}$\pi$ \unit{m} as a function of the free-space wavenumber $k_0$.}
\label{fig:wing_rank}
\end{figure}

In this section, we first aim at assessing the efficiency of the fast direct solver. The rank of the skeleton form $\mat U \mat V^\T$ is the key parameter to observe, as it directly determines the computational complexity of the method, affecting both time and memory required. The first geometry analyzed is the ellipse. Figure \ref{fig:ellipse_rank} shows the rank of the skeleton of the operator, for both TE and TM formulations, evaluated over an ellipse with aspect ratio \num{1.5} and perimeter \num{2}$\pi$ \unit{m}. Secondarily, we have considered an airfoil geometry, resulting from the application of the Joukowsky conformal mapping from the circle, with perimeter \num{2}$\pi$ \unit{m} (\cref{fig:wing_rank}). In both cases, we observe that the rank grows less than linearly with frequency and tends to stabilize to the expected behaviour of $k_0^{1/3}$ in the high frequency limit.  
Consistently, the compression time (i.e., the time required for the skeleton evaluation), dominating the overall inversion time, scales quasi-linearly, as shown in \cref{tab:time}.

\begin{figure}
\centerline{\includegraphics[width=1\columnwidth]{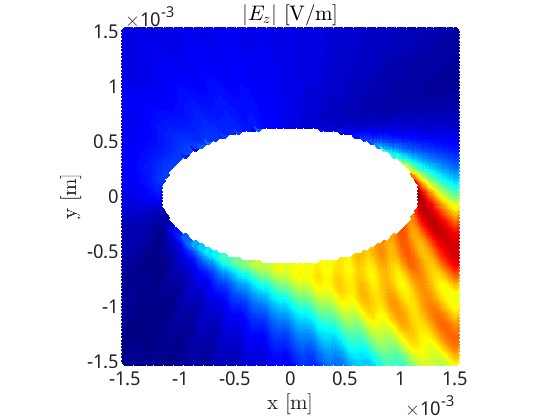}}
\caption{Magnitude of the longitudinal electric field scattered by the skin given a TM plane-wave excitation at the frequency of \SI{1}{THz} impinging at an angle of \num{-1}/\num{5}$\pi$.}
\label{fig:skinfield}
\end{figure}

Then, we applied the solver to the evaluation of the electromagnetic scattering from a skin sample (\cref{fig:skinfield}). In particular, we considered an ellipse of perimeter approximately of \SI{5.85}{mm} excited by a time-harmonic field at the frequency of \SI{1}{THz}. We employed the double Debye model (\cref{eqn:debye}) to approximate the permittivity of the skin, corresponding to a penetration length of approximately \SI{62}{\micro\metre}.

\begin{table}
\vspace{0.2cm}
\centering
\begin{tabular}{|c|c|c|c|}
\hline
$N$     & $k_0$ {[}rad/m{]} & compr. time {[}s{]} & eri \\ \hline
15011 & 833.91        & 59.33            & -      \\ 
30021 & 1667.8         & 129.6             & 1.13 \\ 
45032 & 2501.7        & 210.2            & 1.19  \\
\hline
\end{tabular}
\caption{Compression time required to evaluate the skeleton form $\mat U \mat V^\T$ evaluated over an ellipse with aspect ratio \num{1.5} and perimeter \num{2}$\pi$ \unit{m} at different frequencies, corresponding to different numbers of unknowns $N$, and experimental rate of increase (eri). The simulations have been performed on \num{25} parallel processes.}
\label{tab:time}
\end{table}

\section{Conclusion}
This paper presented a fast direct solver strategy for full-wave modeling and dosimetry at terahertz frequencies. 
This has been obtained by leveraging a preconditioned version of the combined field integral equation, free of spurious high-frequency resonances, and a suitable compression technique for its elliptic spectrum, resulting in an operator matrix invertible in quasi-linear complexity.
The direct nature of the solver makes its use convenient to solve multiple sources problems, where the scatterer response to many different exposures should be analyzed, as it can be the case in dosimetry studies. 

\section*{Acknowledgment} 
The work of this paper has received funding from the European Research Council (ERC) under the European Union’s Horizon 2020 research and innovation programme (grant agreement No 724846, project 321), from the Horizon Europe Research and innovation programme under the EIC Pathfinder grant agreement n° 101046748 (project CEREBRO), and from the ANR Labex CominLabs under the project ``CYCLE''.

\bibliographystyle{IEEEtran}

\bibliography{private_vgiunzioni.bib}

\end{document}